\documentclass{article}%
\usepackage{amsmath}%
\usepackage{amsfonts}%
\usepackage{amssymb}%
\usepackage{graphicx}
\providecommand{\U}[1]{\protect\rule{.1in}{.1in}}

\begin{document}

\title{Longest common subsequences in binary sequences}
\author{John D. Dixon\\Carleton University, Ottawa, Canada}
\maketitle

Given two $\left\{  0,1\right\}  $-sequences $X$ and $Y$ of lengths, say, $m$
and $n$, respectively, we write $L(X,Y)$ to denote the length of the longest
common subsequence (LCS). For example, if $X=01101110$ and $Y=101001011$, then
a longest common subsequence is given by $110111$ ($010111$ is another one),
so $L(X,Y)=6$. We write $X_{k}$ (respectively, $Y_{k}$) to denote the initial
segments of length $k$ for these sequences and for fixed $X$ and $Y$ we shall
write $l(i,j)$ as an abbreviation for $L(X_{i},Y_{j})$ when $i=0,...,m$ and
$j=0,...,n$. It is easily seen that
\begin{equation}
l(i,j)=\left\{
\begin{array}
[c]{l}%
~l(i-1,j-1)+1~\text{if }X_{i}\text{ and }Y_{j}\text{ end in the same symbol}\\
\max\left(  l(i,j-1),l(i-1,j)\right)  \text{ if }X_{i}\text{ and }Y_{j}\text{
end in different symbols}%
\end{array}
\right. \label{def}%
\end{equation}
for $1\leq i\leq m$ and $1\leq j\leq n$. Together with the boundary conditions
$l(i,0)=l(0,j)=0$ for all $i,j$, this gives an efficient way to compute
$L(X,Y)$ for any given $X,Y$ and is generally the way in which $L(X,Y)$ is computed.

\section{Distribution of $L(X,Y)$ for random sequences of same
length\label{sec: diagcase}}

Now suppose that $X$ and $Y$ are random binary sequences of lengths $m$ and
$n$, respectively, and consider the random variable $L(X,Y)$. Let $L(m,n)$
denote its mean value. Then the following are know in the case $m=n$ (see
\cite[Chap. 1]{Stee97}).

\begin{enumerate}
\item The ratio $\gamma_{n}:=L(n,n)/n$ converges to a limit $\gamma$ as
$n\rightarrow\infty$ (see \cite{Chvat75}). The constant $\gamma$ is known as
the Chv\'{a}tal-Sankoff constant, and determination of its value is a
longstanding open problem. The best bounds which have been proved so far are
$0.788071\leq\gamma\leq0.826280$ (see \cite{Luek09}).

\item Numerical evidence suggests that the sequence $\left\{  \gamma
_{n}\right\}  $ is monotonic increasing, but this has not been proved.
However, it is clear that
\[
L(nk,nk)\geq kL(n,n)\text{ for }k=1,2,...
\]
and so $\gamma\geq\gamma_{n}$ for all $n$. As we shall see below, computations
indicate that $\gamma_{16384}$ is approximately $0.81110$ so it seems very
likely that $\gamma>0.81$.

\item The Azuma-Hoeffding inequality shows that
\[
\Pr\left\{  \left\vert L(X,Y)-L(n,n)\right\vert >\lambda\sqrt{n}\right\}
\leq2e^{-\lambda^{2}/8}\text{ for each }\lambda>0
\]
so the values of the random variable $L(X,Y)$ are concentrated around the mean
$L(n,n)$.
\end{enumerate}

\section{Embedding $X$ in a random binary sequence $Y$}

Let $X$ be a fixed binary sequence of length $m$ and consider the probability
$p(X,n)$ that $X$ can be embedded into a random binary sequence of length $n$
for some fixed $n\geq m$, that is, that $L(X,Y)=m$. Since it is equally likely
that $X$ and $Y$ end in the same or different symbols, (\ref{def}) shows that:%
\[
p(X,n)=\frac{1}{2}p(X_{m-1},n-1)+\frac{1}{2}p(X,n-1).
\]
Induction now shows that $p(X,n)$ is independent of the particular sequence
$X$ and depends only on its length, so we can put $p(m,n)$ in place of $p(X,n)
$. In particular, we may assume that $X$ is the sequence of all $1^{\prime}$s
and so $p(m,n)$ is the probability that a random binary sequence of length $n$
has at least $m$ $1^{\prime}$s$.$ Thus, for all $n\geq m$ we have%
\[
p(m,n)=2^{-n}\sum_{k=m}^{n}\binom{n}{k}\text{.}
\]
In particular, the Azuma-Hoeffding inequality shows that%
\[
p(m,n)>1-e^{-\lambda^{2}/8}\text{ if }m<\frac{1}{2}n-\lambda\sqrt{n}
\]
and
\[
p(m,n)<e^{-\lambda^{2}/8}\text{ if }m>\frac{1}{2}n+\lambda\sqrt{n}\text{.}
\]

\section{Distribution of $L(X,Y)$ for sequences of different lengths}

To simplify notation we write $L(r,s)$ (for any two positive reals $r,s$) to
mean $L(\left\lfloor r\right\rfloor ,\left\lfloor s\right\rfloor )$ where
$\left\lfloor ~\right\rfloor $ represents the floor function. Fix $\alpha>0$
and consider the sequence $L(\alpha n,n)/n$ ($n=1,2,...$). With an argument
similar to the argument used to prove the existence of the limit for
$L(n,n)/n$ (see Section \ref{sec: diagcase}) we can show that $L(\alpha n,n)/n
$ converges as $n\rightarrow\infty$ for all $\alpha\geq0$ and denote the limit
by $\psi(\alpha)$. The function $\psi$ has the following properties.

\begin{enumerate}
\item Since $L(\alpha n,n)/n=\alpha L(n,\alpha n)/\alpha n=\alpha
L(\alpha^{-1}r,r)/r$ where $r=\alpha n$, we conclude that $\psi(\alpha
)=\alpha\psi(\alpha^{-1})$ for $\alpha>0$. Clearly $\psi$ is increasing.

\item Since $L(m_{1},n_{1})+L(m_{2},n_{2})\leq L(m_{1}+m_{2},n_{1}+n_{2})$ it
follows that $\lambda\psi(\alpha)+(1-\lambda)\psi(\beta)\leq\psi(\lambda
\alpha+(1-\lambda)\beta)$ when $\alpha,\beta,\lambda$ and $1-\lambda$ are all
nonnegative. Thus $\psi$ is concave (see also \cite{Amsa-etal08}). Since
$\psi$ is bounded and concave in the open interval $(0,\infty)$ it has the
following properties (see \cite{Hard59}): (i) $\psi$ is continuous; (ii)
$\psi$ has a right-hand derivative and a left-hand derivative at each point
with the right-hand derivative not less than the left-hand derivative; and
(iii) these one-sided derivatives are montonic decreasing.

\item It follows from the previous section that $\psi(\alpha)=\alpha$ for
$0\leq\alpha<\frac{1}{2}$ and so by $1.$ we have $\psi(\alpha)=1$ for
$\alpha>2$. It seems that $\psi$ is at least twice differentiable except
perhaps at $\alpha=\frac{1}{2}$ and $\alpha=2$ (see the graph shown below)$.$

\item Let $X$ and $Y$ be two infinite random sequences. Then using the
Azuma-Hoeffding inequality we have for each $\varepsilon>0$:%
\[
\Pr\left\{  \max_{m}\left\vert \frac{L(X_{m},Y_{n})}{n}-\psi(\frac{m}%
{n})\right\vert >\varepsilon\right\}  \rightarrow0\text{ as }n\rightarrow
\infty\text{.}
\]

\item The function $L^{\ast}(m,n):=\sqrt{(4mn-n^{2}-m^{2})/3}$ appears to be a
close approximation to $L(X_{m},Y_{n})$ for $\frac{1}{2}n\leq m\leq2n$ for
\textquotedblleft random\textquotedblright\ $X$ and $Y$. The graph
$m\longmapsto L^{\ast}(m,n)$ is the arc of an ellipse tangential to the line
$y=\frac{1}{2}x$ at $x=n$ and to $y=n$ at $x=2n.$ Its value at $m=n$ is
$L^{\ast}(n,n)=n\sqrt{2/3}$ which is approximately $0.816496n$ and within all
known bounds for $\gamma n$ (\cite{Bund01} claims that $\gamma=0.812653$ but I
suspect that the latter estimate is unreliable). It seems possible that
$\psi(\alpha)$ is equal to $\psi^{\ast}(\alpha):=\sqrt{(4x-x^{2}-1)/2}$ and
$\gamma=\sqrt{2/3}$. I conjecture that at any rate $\psi^{\ast}$ is an upper
bound to $\psi$.
\end{enumerate}

Computing $L(X,Y)$ for random $X,Y$ we obtained the graph in Figure 1.%

\begin{figure}[ptb]%
\centering
\includegraphics[
height=4.4691in,
width=5.9493in
]%
{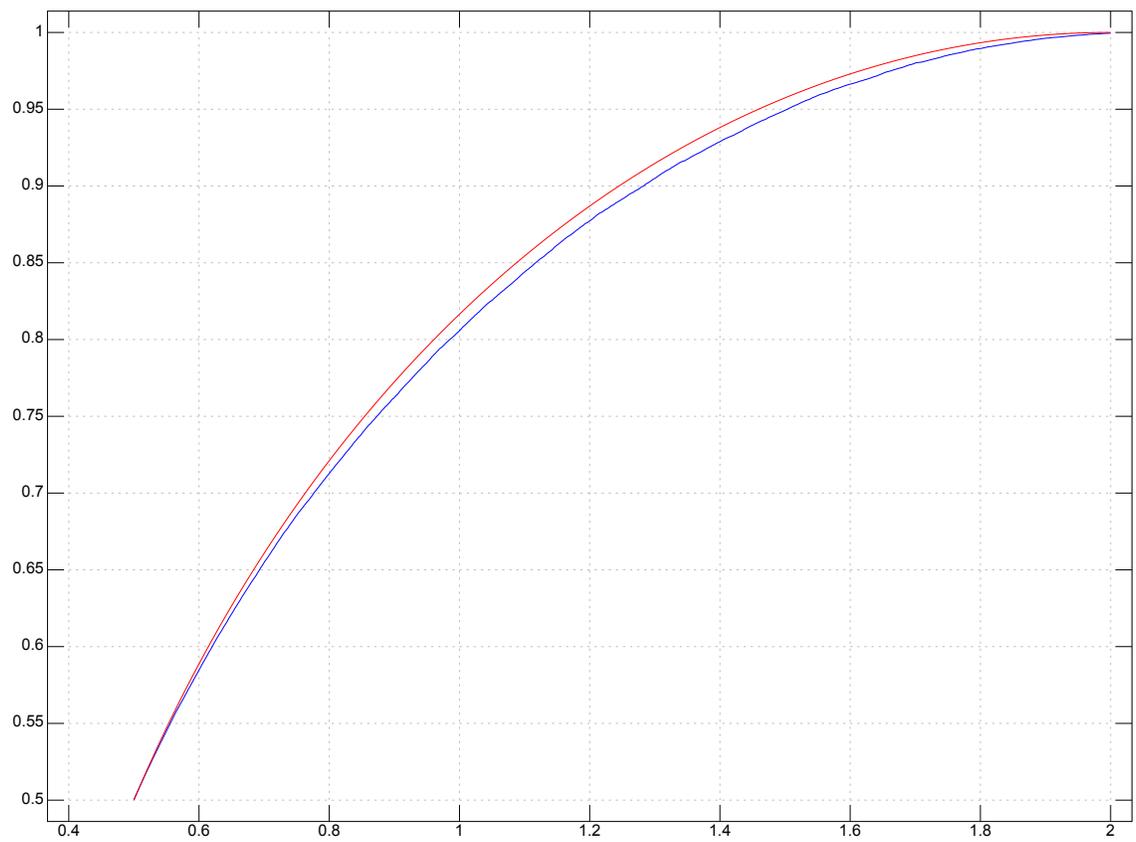}%
\caption{$L^{\ast}(m,1000)/1000${\protect\small \ and }$L(m,1000)/1000$%
{\protect\small \ with }$0.5<m/1000<2$}%
\end{figure}

\qquad

\section{From the discrete to the infinite\label{sec: disinf}}

If $X$ and $Y$ are infinite random binary sequences, then (\ref{def}) can be
used to compute the values of $l(m,n)$ ($=L(X_{m},Y_{n})$). Let $l[m]$ denote
the $m$th row of the infinite matrix $\left[  l(m,n)\right]  _{k,n=0,1,...}$ .
For any infinite (real) vector $v=(v(i))_{i=0,1,....}$ we define the
\emph{maximizer }$%
{\displaystyle\bigvee}
$ such that $%
{\displaystyle\bigvee}
v$ is the vector whose $j$th entry is the maximum of the $v(i)$ for $i\leq j$.

Define $T$ and $\bar{T}$ as operators on vectors by%
\[
(Tv)(i):=\left\{
\begin{array}
[c]{c}%
v(i-1)+1\text{ if }y_{i}=1\\
v(i)\text{ otherwise}%
\end{array}
\right.
\]
where $y_{i}$ is the $i$th entry of $Y$ and%
\[
(\bar{T}v)(i):=\left\{
\begin{array}
[c]{c}%
v(i-1)+1\text{ if }y_{i}=0\\
v(i)\text{ otherwise}%
\end{array}
\right.  .
\]
It follows from (\ref{def}) that
\[
l[m]=%
{\displaystyle\bigvee}
Tl[m-1]\text{ if }x_{m}=1
\]
and
\[
l[m]=%
{\displaystyle\bigvee}
\bar{T}l[m-1]\text{ if }x_{m}=0.
\]
This defines the successive rows of the LCS table using global operations. In
particular, $l[m]$ is obtained from $l[0]$ by applying the nonlinear operators
$%
{\displaystyle\bigvee}
T$ and $%
{\displaystyle\bigvee}
\bar{T}$ in random order. In view of the relationship between the values of
$l[m]$ and the values of $\psi$, this may give a hint as to the kinds of
operators which leave $\psi$ fixed, and perhaps $\psi$ may be determined in
this way.

Computations for this paper were done using the $J$-language developed by
Iverson and Hui (see \cite{J}). Although $J$ is an interpreted language it is
fast because it is based on a large number of carefully integrated and
optimized subroutines. The most efficient programs in $J$ turned out to be
based on the global approach described above. $J$ is a very concise language
and the full program to compute $L(X_{m},Y)$ ($m=1,2,...$) for two finite
$\left\{  0,1\right\}  $-lists $X$ and $Y$ is given as follows (lines
beginning $NB.$ are comments):
\begin{verbatim}
LCS=: 3 : 0 "_ 0 _
NB. LCS Y; X  returns a list of the lengths of the LCS 
NB. for Y and all initial segments of X
'u v' =. y
val=. v*0
for_e. u do.
 val=. >./\(e=v)} val,: }: 1, 1+val
end.
val
)
 
\end{verbatim}

We estimated $L(n,n)$ by taking the mean value of $50$ trials of
$\mathtt{LCS}$ $\mathtt{Y;X}$\texttt{\ where }$X$ and $Y$ were random
$\{0,1\}$-lists of length $n$ ($err$ is the standard deviation for the sample mean):

\begin{flushleft}%
\begin{tabular}
[c]{llllllllll}%
\emph{n} & {\scriptsize 64} & {\scriptsize 128} & {\scriptsize 256} &
{\scriptsize 512} & {\scriptsize 1024} & {\scriptsize 2048} &
{\scriptsize 4096} & {\scriptsize 8192} & {\scriptsize 16384}\\
\emph{L(n,n)} & {\scriptsize 0.77406} & {\scriptsize 0.78266} &
{\scriptsize 0.79656} & {\scriptsize 0.80121} & {\scriptsize 0.80594} &
{\scriptsize 0.80711} & {\scriptsize 0.80942} & {\scriptsize 0.81031} &
{\scriptsize 0.81110}\\
\emph{err} & {\scriptsize 0.00467} & {\scriptsize 0.00286} &
{\scriptsize 0.00166} & {\scriptsize 0.00108} & {\scriptsize 0.00061} &
{\scriptsize 0.00052} & {\scriptsize 0.00032} & {\scriptsize 0.00021} &
{\scriptsize 0.00014}%
\end{tabular}

\end{flushleft}

\section{Generating LCS table with a finite state machine}

For any two $(0,1)$-sequences $X$ and $Y$ we consider the (possibly infinite)
table whose $(i,j)$th entry is $L(X_{i},Y_{j})$ ($i,j=0,1,...$). The initial
row and column of this table consists of $0$'s and we can define a table $T$
with $(0,1)$-entries by
\[
T_{i,j}:=L(X_{i},Y_{j})-L(X_{i-1},Y_{j})\text{ for }i=1,2,...\text{ and
}j=0,1,...\text{ .}
\]
Evidently knowledge of the entries of $T$ determine the values of
$L(X_{i},Y_{j})$. We can compute the rows of $T$ recursively with a finite
state machine as follows.

To compute values of $T_{ij}$ with given $j>0$ and $i=1,2,...$ we use the
triple $(T_{i-1,j-1}$, $T_{i,j-1},f)$ where $f$ is a flag equal to $0$ or $1$
which defines the state of the machine. As input we have the pair $T_{i-1,j}$
and $Y_{j}$. The machine computes $T_{ij}$, moves into a new state defined by
$(T_{i-1,j},$ $T_{i,j},\tilde{f}$ $)$ and outputs the value of $T_{ij}$. The
flag represents the carry which is necessary when the maximizer is applied to
the row in computation of $L(X_{i},Y_{j})$ described in Section
\ref{sec: disinf}. In this form the finite state machine requires $2^{3}$
states, but some of these turn out to be indistinguishable so we can reduce to
four states. We do not give the details but provide the final tables for a fsm
(see Tables 1 and 2).

\begin{center}%
\begin{table}[tbp] \centering
\begin{tabular}
[c]{c|cccc}
& 00 & 01 & 10 & 11\\\hline
0 & 0 & 0 & 1 & 1\\
1 & 1 & 0 & 1 & 1\\
2 & 0 & 2 & 1 & 3\\
3 & 1 & 2 & 1 & 3
\end{tabular}
\caption
{State Transition Table.  Rows are labelled by states and columns by pairs $(T_{ij}, Y_j)$.}%
\end{table}%

\bigskip%

\begin{table}[tbp] \centering
\begin{tabular}
[c]{c|cccc}
& 00 & 01 & 10 & 11\\\hline
0 & 0 & 0 & 0 & 0\\
1 & 0 & 1 & 1 & 1\\
2 & 0 & 0 & 0 & 0\\
3 & 0 & 1 & 1 & 1
\end{tabular}
\caption{Output Table}%
\end{table}%

\end{center}

\section{Partially ordered sets and longest chains}

Another way to describe the same problem is as follows. Given two binary
sequences $X$ and $Y$ of lengths $m$ and $n$, respectively, we define the set
$P:=\left\{  (i,j)~|~x_{i}=y_{j}\text{ with }1\leq i\leq m\text{ and }1\leq
j\leq n\right\}  $. We partially order $P$ with $<$ where $(i,j)<(i^{\prime
},j^{\prime})$ $\iff$ $\left[  i<i^{\prime}\text{ and }j<j^{\prime}\right]  $.
It can be verified that $(i_{1},j_{1})<(i_{2},j_{2})<\ldots<(i_{k},j_{k})$ is
a chain in $(P,<)$ if and only if $(x_{i_{1}},x_{i_{2}},...,x_{i_{k}%
})=(y_{j_{1}},y_{j_{2}},...,y_{j_{k}})$ is a common subsequence of $X$ and
$Y$. In particular, the longest chain in $(P,<)$ has length $L(X,Y)$.

\end{document}